\documentclass{article}
\usepackage[margin=1in]{geometry}
\usepackage[utf8]{inputenc}
\usepackage{natbib}
\usepackage{amsmath}
\usepackage{mathtools}
\usepackage{amsthm}
\usepackage{amsfonts}
\usepackage{amssymb}
\usepackage{graphicx} 
\usepackage{textcomp}
\usepackage{hyperref}
\usepackage{color}
\usepackage[dvipsnames]{xcolor}
\usepackage{enumerate}
\usepackage[inline]{enumitem}
\usepackage[explicit]{titlesec}
\usepackage{xifthen}
\usepackage{xspace}
\usepackage{comment}

\bibpunct{[}{]}{,}{n}{}{,}

\definecolor{LinkColor}{HTML}{5c5c84}
\definecolor{CiteColor}{HTML}{5c5c84}
\definecolor{URLColor}{rgb}{0.3,0.5,0.7}
\hypersetup{
    colorlinks=true,
    linkcolor=LinkColor,
    filecolor=magenta,
    urlcolor=URLColor, 
    citecolor=CiteColor,
    anchorcolor=blue
}

\allowdisplaybreaks

\DeclareSymbolFont{bbold}{U}{bbold}{m}{n}
\DeclareSymbolFontAlphabet{\mathbbold}{bbold}

\def\cmtoutcolorvalue{0.7}
\definecolor{cmtoutcolor}{gray}{\cmtoutcolorvalue}

\newcommand{\Enum}[2][]{\ifthenelse{\isempty{#1}}{(#2)}{\Label{#1}{(#2)}}\xspace}
\newcommand{\TagEqn}{\stepcounter{equation}\tag{\theequation}}

\MakeRobust{\ref}

\makeatletter
\newcommand{\Label}[2]{%
  \@bsphack
  \csname phantomsection\endcsname 
  #2%
  \def\@currentlabel{#2}{#1}
  \@esphack
}
\makeatother

\def\ldelim{(}
\def\rdelim{)}

\definecolor{todocolor}{rgb}{1,0,0}
\newcommand{\xToDo}[1]{\colorbox{yellow}{\normalfont%
 $\ldelim$\textbf{To\,do}\ifthenelse{\isempty{#1}}{\textbf{:}\;Fill this in.}{\textbf{:}\;{#1}}$\rdelim$}}
\newcommand{\ToDo}[1]
{%
  \colorlet{origcolor}{.}
  \color{todocolor}%
  \relax\ifmmode\text{\xToDo{#1}}\else\xToDo{#1}\fi
  \color{origcolor}%
}%

\newcommand{\nm}[1]{{\color{black}#1}}

\makeatletter
\titleformat{\paragraph}[runin]{\normalfont\bfseries}{}{}{#1. }
\makeatother

\makeatletter
\renewenvironment{proof}[1]{\par
  \pushQED{\qed}%
  \normalfont \topsep6\p@\@plus6\p@\relax
  \trivlist
  \item\relax
  {\itshape
  \proofname\ifthenelse{\equal{#1}{}}{\@addpunct{.}}{ of #1\@addpunct{.}}}\hspace\labelsep\ignorespaces
}{%
  \popQED\endtrivlist\@endpefalse \vskip8\p@\@plus6\p@\relax
}
\makeatother

\makeatletter
\newenvironment{subproof}[1]{\par
  \pushQED{\qed}
  \normalfont \topsep6\p@\@plus6\p@\relax
  \trivlist
  \item\relax
  {\itshape
  \proofname\ifthenelse{\equal{#1}{}}{\@addpunct{.}}{ of #1\@addpunct{.}}}\hspace\labelsep\ignorespaces
}{%
  \popQED\endtrivlist\@endpefalse \vskip8\p@\@plus6\p@\relax
}
\makeatother



\newtheoremstyle{claim}
{3pt}
{3pt}
{\itshape}
{}
{\bfseries}
{.}
{.5em}
{}

\newtheoremstyle{property}
{3pt}
{3pt}
{\itshape}
{}
{}
{}
{.5em}
{}

\newtheoremstyle{LemmaRepeat}
{3pt}
{3pt}
{\itshape}
{}
{\bfseries}
{.}
{.5em}
{\thmname{#1}\thmnote{ #3}}

\newtheoremstyle{mystyle}
  {}
  {}
  {\itshape}
  {}
  {\bfseries}
  {.}
  { }
  {\thmname{#1}\thmnumber{ #2}\thmnote{ (#3)}}

\theoremstyle{mystyle}

\newtheorem{theorem}{Theorem}[section]

\newtheorem{lemma}[theorem]{Lemma}
\newtheorem{prop}[theorem]{Proposition}

\newtheorem{remark}{Remark}[section]
\newtheorem*{theorem*}{Restatement of Theorem}
\newtheorem*{lemma*}{Restatement of Lemma}
\newtheorem*{claim*}{Restatement of Claim}
\newtheorem*{corollary*}{Restatement of Corollary}
\newtheorem*{fact*}{Restatement of Fact}
\newtheorem{definition}{Definition}[section]
\theoremstyle{claim}
\newtheorem{claim}[theorem]{Claim}
\theoremstyle{property}

\theoremstyle{LemmaRepeat}



\newcommand{\Theorem}{Theorem~}
\newcommand{\Lemma}{Lemma~}

\newcommand{\Claim}{Claim~}

\newcommand{\Equation}{Equation~}
\newcommand{\Section}{Section~}

\newcommand{\Step}{Stage~}

\newcommand{\Definition}{Definition~}

\newcommand{\Lemmas}{Lemmas~}

\newcommand{\Equations}{Equations~}
\newcommand{\Sections}{Sections~}

\newcommand{\Steps}{Stages~}

\newcommand{\etal}[1][~]{et al.#1}

\newcommand{\Th}{^{\text{th}}}

\newcommand{\AlnSp}{\ensuremath{~~~~}}
\newcommand{\NegAlnSp}{\negphantom{\AlnSp}}

\newcommand{\negphantom}[1]{\ifmmode\settowidth{\dimen0}{$#1$}\else\settowidth{\dimen0}{#1}\fi\hspace*{-\dimen0}}

\newcommand{\R}{\mathbb{R}}
\newcommand{\Z}{\mathbb{Z}}
\newcommand{\F}{\mathbb{F}}
\newcommand{\Nat}{\mathbb{N}}

\renewcommand{\Vec}[1]{\boldsymbol{\mathrm{#1}}}

\renewcommand{\Pr}{P}
\newcommand{\E}{\operatornamewithlimits{\mathbb{E}}}

\newcommand{\Uniform}{\operatorname{Uniform}}

\newcommand{\sign}{\operatorname{sign}}

\newcommand{\I}{\mathbb{I}}

\newcommand{\defeq}{\triangleq}

\newcommand{\T}{T}

\newcommand{\Mat}[1]{\mathbf{#1}}

\newcommand{\nequiv}{\not\equiv}

\newcommand{\n}{n}
\renewcommand{\t}{t}
\newcommand{\m}{m}
\newcommand{\f}{f}
\newcommand{\fx}[1]{f_{#1}}
\newcommand{\x}{x}
\newcommand{\ux}{u}

\newcommand{\cconst}{C}

\newcommand{\mlb}{\cconst \t}

\newcommand{\block}{B}
\newcommand{\blockset}{\mathbb{B}}
\newcommand{\Xset}{X}
\newcommand{\w}{w}

\newcommand{\Wset}{W}
\newcommand{\Vset}{V}

\newcommand{\Rset}{R}
\newcommand{\Sset}{S}
\newcommand{\Tset}{T}
\newcommand{\Uset}{U}
\newcommand{\Yset}{Y}

\newcommand{\pset}{\mathcal{P}}
\newcommand{\p}{p}
\newcommand{\minp}{\frac{\t}{4}}
\newcommand{\maxp}{\frac{5\t}{12}}

\newcommand{\ix}{i}
\newcommand{\jx}{j}
\newcommand{\kx}{k}

\newcommand{\lx}{\ell}

\newcommand{\Ix}{I}

\newcommand{\wtop}[1]{\phi_{#1}}

\renewcommand{\k}{k}
\newcommand{\h}{h}
\newcommand{\Uinterval}[1]{U_{#1}}

\newcommand{\Vyinterval}[2]{V_{#2}^{#1}}

\newcommand{\yx}{y}

\newcommand{\rv}{r}
\newcommand{\sv}{s}
\newcommand{\uv}{u}
\newcommand{\vv}{v}
\newcommand{\crsuv}[4]{c'_{#1,#2,#3,#4}}

\newcommand{\crsu}[3]{c_{#1,#2,#3}}

\newcommand{\MODPVALUES}{\frac{\t}{2}+\yx\p : \yx \in \YVALUESSET}
\newcommand{\YVALUESSET}{\{ -1,0,1,2 \}}
\newcommand{\NYVALUES}{4}

\newcommand{\jO}{j_{\ast}}
\newcommand{\aconst}{a}
\newcommand{\aconstvalue}{8}
\newcommand{\bconst}{b}
\newcommand{\bconstvalue}{60}
\newcommand{\alphaconst}{c}
\newcommand{\alphaconstvalue}{\frac{1}{96}}

\newcommand{\TMARGIN}{\t}
\newcommand{\jOset}{J}
\newcommand{\jxx}{j'}
\newcommand{\Iyjl}[3]{I_{#2,#3}^{#1}}
\newcommand{\IyjlX}[2]{\tilde{I}_{#2}^{#1}}
\newcommand{\Zyj}[2]{Z_{#2}}
\newcommand{\Zyjl}[3]{Z_{#2,#3}^{#1}}
\newcommand{\xv}{x}
\newcommand{\yv}{y}
\newcommand{\tO}{t_{0}}
\newcommand{\Aset}{A}
\newcommand{\Bset}{B}
\newcommand{\Zset}{Z}

\newcommand{\Vleqt}{V_{\leq \t}}
\title{Lower bound on the size of splitting systems}
\author{Fan Chung \and Namiko Matsumoto \and Arya Mazumdar\thanks{The authors are with University of California San Diego. This work is supported in part by NSF award 2217058 and NSF GRFP
No. DGE-2038238.}}
\date{\today}
\begin{document}
\maketitle
\begin{abstract}
For even $\t$ and an $\n$-set $\Xset$, an $(\n,\t)$-splitting system is a family ${\blockset} \subseteq 2^\Xset$ with the property that for each subset $\Yset \in {\Xset \choose \t}$, there exists $\block \in \blockset$ such that $|\block \cap \Yset| = \t/2$. The problem of interest is to determine the minimum size of such splitting systems. We show that
any $(\n,\t)$-splitting system has size at least of order $\t$, improving the previous lower bound of $\sqrt \t$ while the best known upper bound is  of order $\t^{3/2}$.

\end{abstract}

\section{Introduction}
\label{section:intro}

Splitting systems, as defined below, were formally introduced by Stinson in a 2002 paper \cite{stinson2002some} on a generalization of Coppersmith algorithm \cite{coppersmith2003minimum,menezes2018handbook} for the low Hamming weight discrete logarithm problem. 

\begin{definition}[splitting system] 
\label{def:} 
Let \(  \m, \n,\t \in \Nat  \), where \(  \t  \) is even.
Let \(  \Xset  \) be a set with \(  |\Xset| = \n  \), and let \(  \blockset \subseteq 2^{\Xset}  \) be a family of \(  \m  \) \emph{blocks} \(  \block_{1}, \dots, \block_{\m} \subseteq \Xset  \).
The \emph{set system} \(  (\Xset, \blockset)  \) is an \emph{\(  (\n,\t)  \)-splitting system} if for every \(  \Yset \in \binom{\Xset}{\t}  \), there exists a block \(  \block_{\ix} \in \blockset  \) such that
\(  |\Yset \cap \block_{\ix}| = \frac{\t}{2}  \). 
\end{definition} 

A natural problem is to ask for the splitting system $(\Xset, \blockset)$ so that the cardinality
of $\blockset$ is minimized.

Several additional works have studied splitting systems
\citep{deng2007constructions,deng2004lovasz,deng2011stein,ling2004splitting,roh2012constructions} and proposed generalizations \citep{condon2014generalizations,coskey2018splitting} for their merit beyond
the original application in reference \cite{stinson2002some}.
Moreover, splitting systems have connections to other combinatorial families and problems, such as separating families, balanced perfect hash families, sunflowers, and balancing sets of binary vectors \citep{alon2020unbalancing,condon2014generalizations,deng2007constructions,frankl1987forbidden,hegedHus2010balancing,ling2004splitting}.
Splitting systems can be regarded as special cases of splitters, which are combinatorial structures \cite{naor1995splitters} that partition every \(  \t  \)-subset into a specified number of equally-sized parts (up to rounding),
often used as a derandomization technique 
\citep{alon2010balanced,alon2006algorithmic,burjons2025uniform,lokshtanov2021approximate,naor1995splitters}.
Explicit constructions of splitters can be used to obtain explicit constructions of families of perfect hash functions \cite{alon2010balanced}.
Nevertheless, the sizes of explicit constructions tend to significantly exceed what can be shown to exist nonconstructively.
\color{black}

Both constructive and nonconstructive methods have been proposed to establish splitting systems using a variety of methods including direct or recursive constructions \citep{deng2007constructions,roh2012constructions,stinson2002some}, the probabilistic method \citep{stinson2002some}, perfect hash families \citep{stinson2002some}, separating systems \citep{ling2004splitting}, the Lov\'{a}sz local lemma \citep{deng2004lovasz}, and the Stein-Lov\'{a}sz theorem \citep{deng2011stein}.
Coppersmith (see \cite[\Lemma 2.1]{stinson2002some}) established existence of splitting systems with \(  \m = \frac{\n}{2}  \) blocks for any even \(  \n  \) and \(  \t  \), while 
Deng \etal \cite{deng2007constructions} showed existence of \(  (\n,3)  \)-splitting systems with \(  \m = \lfloor \frac{\n}{2} \rfloor  \), where the definition of splitting systems extends to odd-valued \(  \t  \) by replacing \(  \frac{\t}{2}  \) with \(  \lfloor \frac{\t}{2} \rfloor  \).
For \(  \t \in \{2,4\}  \) and for \(  \t = 3  \), Ling \etal \cite{ling2004splitting} and, respectively, Roh and Hahn \cite{roh2012constructions} proved that splitting systems with \(  \m \asymp \log \n  \) blocks exist.
However,
for arbitrary \(  \t > 4  \), only splitting systems with significantly higher dependence of \(  \m  \) on \(  \t  \) are known to exist.
Stinson \cite{stinson2002some} presented an explicit, recursive construction based on balanced perfect hash families with \(  \m \lesssim \t^{4} \log^{2} \n  \) blocks,
and later, Deng \etal \cite{deng2011stein} obtained a constructive upper bound of
{\(  \m \lesssim \t^{3/2} \log \n  \)}
through the Stein-Lov\'{a}sz algorithm.
Various nonconstructive methods have also shown an
\(  \m \lesssim \t^{3/2} \log \n  \)
bound \citep{deng2004lovasz,stinson2002some}, which is the best upper bound on the minimal size of splitting systems obtained so far for \(  \t \gg 1  \).

For lower bounds on the size  $\m$ of splitting systems, it is shown that  $ \m \gtrsim \log \n $ blocks are necessary for $\t \in \{2,4\}$  (see \cite{ling2004splitting}) and for general $t$ (see \cite{deng2007constructions}).
A naive counting bound
gives a lower bound of \(  \m \gtrsim \sqrt{\t}  \) in the regime \(  \t \gg 1  \).
This follows from standard facts (see \Lemma \ref{lemma:binom-coef-ineq:1}) that imply that a single block can split at most \(  \frac{1}{\sqrt{\t}}  \)-fraction of the \(  \t  \)-subsets of \(  \Xset  \).
This leaves an order-\(  \t  \) separation for the dependence on \(  \t  \) between the known lower and upper bounds in the asymptotic regime.
This note contributes a substantial reduction to only a
\(  \sqrt{\t}  \)-gap
by showing that \(  \m \gtrsim \t  \) blocks are necessary for an \(  (\n,\t) \)-splitting system for all \(  \n  \) and all \(  \t  \) exceeding a certain constant when \(  \t \leq \frac{\n}{2}  \).

It is worth mentioning that Coskey \etal \cite{coskey2018splitting} obtained a tight bound of \(  \m \asymp \n  \), but that work proposed a different definition for splitting systems,
{and their lower bound does not seem to immediately extend to the \(  (\n, \t)  \)-splitting systems that are the subject of this note.}
{%
The same-order tight bound, \(  \m \asymp \n  \), has been shown several times in the literature for \(  (\n, \frac{\n}{2})  \)-splitting systems \citep{alon2020unbalancing,enomoto1987codes,frankl1987forbidden,hegedHus2010balancing}.
Notably, while working on this note, we observed that for odd \(  \frac{\t}{2}  \), the proof in reference \cite{enomoto1987codes} can be combined with a lemma in reference \cite{bassalygo2000codes} to prove the same-order lower bound as obtained here, \(  \m \gtrsim \t  \), with a more streamlined argument based on linear codes with given distances and with a substantially smaller constant in the bound.
However, the argument breaks down when \(  \frac{\t}{2}  \) is even, whereas this paper requires no assumption on the parity of \(  \frac{\t}{2}  \).
}

\subsection{Connections to 1-bit compressed sensing and threshold group testing}


The original motivation for splitting systems was the development of algorithms to solve a cryptographic problem \citep{stinson2002some}.
However, problems in statistics and signal processing have a close connection to splitting systems, such as exact recovery of sparse binary vectors in 1-bit compressed sensing, and  threshold group testing.

In the universal exact recovery of binary vectors problem in 1-bit compressed sensing, every \(  \t  \)-sparse binary vector---a vector \(  \Vec{x} \in \{0,1\}^{\n}  \) with at most \(  \t \ll \n  \) nonzero entries---is exactly reconstructed from
\(  \m  \) one-bit quantized linear projections
\(  \Vec{y} = \sign( \Mat{A} \Vec{x} ) \in \{-1,1\}^{\m}  \), where \(  \Mat{A} \in \R^{\m \times n}  \) is an  \(  \m \times \n  \) real-valued \emph{measurement matrix}
whose rows form \(  \m  \) \emph{measurements}, and the \(  \sign  \) function maps entries to their \(  \pm  \)-signs coordinate-wise.
The folklore lower bound on the number \(  \m  \) of measurements necessary in this problem is \(  \m \geq \log_{2} \binom{\n}{\t} \gtrsim \t \log \frac{\n}{\t}  \), which follows from the simple information-theoretic argument that there are \(  2^{\m}  \)-many \(  \m  \)-dimensional binary vectors and \(  \binom{\n}{\t}  \)-many binary vectors with Hamming weight \(  t  \).
However, despite some efforts \citep{bansal2022universal,flodin2019superset}, to date, the best upper bound for this problem, achieved probabilistically, is \(  \m \lesssim \t^{3/2} \log \frac{\n}{\t}  \) \citep{bansal2022universal,flodin2019superset,jacques2013robust,matsumoto2023improved}, which exceeds this lower bound by a factor of order-\(  \sqrt{\t}  \).
Works such as references \cite{bansal2022universal} and \cite{flodin2019superset} have noted the difficulty of this problem: indeed, standard techniques for measurement designs and algorithms that have succeeded in other problems in 1-bit compressed sensing
fail to
rectify
this discrepancy between the lower and upper bounds.
In fact, the same \(  \sqrt{\t}  \)-gap
persists even under the restriction to arbitrary \(  \{-1,1\}  \)-valued measurement designs, where in particular, Bansal \etal \cite{bansal2022universal} obtained a probabilistic upper bound of \(  \m \lesssim \t ^{3/2}\log \frac{\n}{\t}  \) by taking the entries in \(  \Mat{A}  \) to be i.i.d. and uniformly in \(  \{-1,1\}  \).

For binary vectors with Hamming weight exactly \(  \t  \), this specialization is equivalent to the universal recovery problem in threshold group testing with a \(  \frac{\t}{2}  \)-threshold.
Universal recovery in threshold group testing is a related problem in which each vector \(  \Vec{x} \in \{0,1\}^{\n}  \) with Hamming weight exactly \(  \t  \) is reconstructed from a vector of \(  \m  \) \(  \{0,1\}  \)-valued test outcomes \(  \Vec{y} = T_{w}(\Mat{A} \Vec{x})  \), where \(  \Mat{A} \in \{0,1\}^{\m \times \n}  \) is a \emph{test matrix}, and where for a fixed threshold \(  w \in [\t-1] \), the threshold function \(  T_{w}  \) is applied entry-wise,
mapping values above the threshold \(  w  \) to \(  1  \) and below \(  w  \) to \(  0  \).
With a threshold of \(  w = \frac{\t}{2}  \), \(  \sign(\Vec{u}) = 2 T_{w}(\Vec{u}) - \Vec{1}  \) for every \(  \Vec{u} \in \R^{\m}  \), giving the equivalence with the special case of the 1-bit sensing problem mentioned above.
For arbitrary thresholds \(  w \in [\t-1]  \), an analogous equivalence arises for \(  \{ -\frac{w}{t-w}, 1 \}  \)-valued measurement matrices in the 1-bit sensing problem. 
In threshold group testing with threshold \(  \frac{\t}{2}  \), the best known lower bound is, again, the information-theoretic bound: \(  \m \gtrsim \t \log \frac{\n}{\t}  \),
while the known upper bound, a probabilistic result, is \(  \m \gtrsim \t^{3/2} \log \frac{\n}{\t}  \) \citep{de2020subquadratic}.
For general thresholds, De Marco \etal \cite{de2020subquadratic} obtained a lower bound of \(  \m \gtrsim \min \{ \frac{\t^{2}}{w^{2}}, \frac{\n}{w} \}  \), where \(  w  \) denotes the threshold. However, for the regime relevant to this note, \(  w = \frac{\t}{2}  \), their result reduces to a lower bound of \(  \m \gtrsim 1  \). 

The \(  \frac{\t}{2}  \)-threshold group testing and, under the specialization to \(  \Mat{A} \in \{-1,1\}^{\m \times n}  \), the 1-bit sensing universal exact recovery problems have a fundamental connection to 
splitting systems,
as illustrated in the following lemma that establishes that the \(  \m  \) sets of column indexes of the \(  1  \)s in each row of \(  \Mat{A}  \) are necessarily in correspondence with an \(  (\n,\t)  \)-splitting system.

\begin{prop} 
\label{prop:1bcs} 
Let \(  \t \in \Nat  \) be even.
For each \(  \ix \in [\m]  \), let \(  \block_{\ix} \subseteq [\n]  \) be the set indexing the \(  1  \)-valued entries in the \(  \ix\Th  \) row of \(  \Mat{A}  \), and write \(  \blockset = \{ \block_{1}, \dots, \block_{\m} \}  \) and \(  \Xset = [\n]  \).
Denote the set of \(  \t  \)-sparse binary vectors
by \(  \Vleqt \defeq \{  \Vec{x} \in \{0,1\}^{\n} : |\Vec{x}|_{H} \leq \t  \}  \),
where \(  |\Vec{x}|_{H}  \) denotes the Hamming weight of a vector \(  \Vec{x}  \). 
Suppose that the map \(  \sign(\Mat{A} \cdot)  \) is an injection over \(  \Vleqt  \).
Then, \(  (\Xset, \blockset)  \) is an \(  (\n,\t)  \)-splitting system.
\end{prop} 

\begin{proof}{Proposition \ref{prop:1bcs}} 
Suppose for the sake of contradiction that a measurement matrix \(  \Mat{A} \in \{-1,1\}^{\m \times \n}  \) has the property that \(  \sign(\Mat{A} \cdot)  \) is an injection over \(  \Vleqt  \),
but \(  (\Xset, \blockset)  \) is not an \(  (\n,\t)  \)-splitting system.
Then, there exists \(  \Yset \in \binom{\Xset}{\t}  \) such that for all \(  \ix \in [\m]  \), \(  ||\Yset \cap \block_{\ix}| - \frac{\t}{2}| \geq 1  \).
Let \(  \Vec{x} \in \Vleqt  \) be a vector whose support (the set indexing its nonzero entries) is exactly \(  \Yset  \), and let \(  \Vec{x'} \in \Vleqt  \) be the vector whose support is exactly \(  \Yset \setminus \{j\}  \) for an arbitrary choice of \(  \jx \in \Yset  \).
Then, for each \(  \ix \in [\m]  \), either
\Enum[\label{case:1bcs:i}]{i} \(  (\Mat{A} \Vec{x})_{\ix} \geq 2  \), or
\Enum[\label{case:1bcs:ii}]{ii} \(  (\Mat{A} \Vec{x})_{\ix} \leq -2  \).
In the case of \ref{case:1bcs:i}, \(  (\Mat{A} \Vec{x'})_{\ix} \geq (\Mat{A} \Vec{x})_{\ix} - 1 \geq 1  \), and hence, \(  (\sign(\Mat{A} \Vec{x}))_{\ix} = (\sign(\Mat{A} \Vec{x'}))_{\ix} = 1  \).
Otherwise, in the case of \ref{case:1bcs:ii}, \(  (\Mat{A} \Vec{x'})_{\ix} \leq (\Mat{A} \Vec{x})_{\ix} + 1 \leq -1  \), and hence, \(  (\sign(\Mat{A} \Vec{x}))_{\ix} = (\sign(\Mat{A} \Vec{x'}))_{\ix} = -1  \).
Therefore, \(  \sign(\Mat{A} \Vec{x}) = \sign(\Mat{A} \Vec{x'})  \) but \(  \Vec{x} \neq \Vec{x'}  \), implying the contradiction that \(  \sign(\Mat{A} \cdot)  \) is not injective.
\end{proof} 


Proposition \ref{prop:1bcs} provides a means to lower bound the number of measurements necessary in these specializations of the 1-bit sensing and threshold group testing problems.
On the other hand, distinguishing pairs of \(  \t  \)-sparse binary vectors with Hamming distance \(  1  \) or \(  2  \) appears to be a bottleneck for solving the recovery problems---which is indicated by the probabilistic upper bounds \citep{bansal2022universal,de2020subquadratic}---suggesting that optimal \(  (\n,\t)  \)-splitting system constructions may additionally guide the development of optimal measurement designs and reconstruction algorithms for these two recovery problems.

\section{Main result}
\label{section:main-thm}

The main theorem of this note is the following lower bound on the minimum number of blocks \(  \m  \) necessary for a set system to be a splitting system.

\begin{theorem} 
\label{thm:main-thm-ss-lb:1} 
Let \(  \m, \n, \t \in \Z_{+}  \) such that \(  \t  \) is even and \(  \tO \leq \t \leq \frac{\n}{2}  \) for an absolute constant \(  \tO > 0  \).
Let \(  \blockset = \{ \block_{1}, \dots, \block_{\m} \subseteq \Xset \}  \), where \(  |\Xset| = \n  \).
Suppose \(  (\Xset, \blockset)  \) forms an \(  (\n,\t)  \)-splitting system.
Then, $\m \geq \cconst \t$ for an absolute constant $\cconst$.
\end{theorem} 

The constant $\cconst$ can be taken as defined in Definition \ref{def:constants}.

Before we proceed to prove \Theorem \ref{thm:main-thm-ss-lb:1}, we give an outline of the proofs of \Theorem \ref{thm:main-thm-ss-lb:1} and the main technical lemma.


\subsection{Overview of the proof of \Theorem \ref{thm:main-thm-ss-lb:1}}
\label{section:pf-overview:main-thm}

The proof of \Theorem \ref{thm:main-thm-ss-lb:1} extends techniques in  references \cite{alon2020unbalancing,hegedHus2010balancing} for $(\n, \frac{\n}{2})$-splitting systems.
The proof relies on two technical lemmas -- formally stated in \Section \ref{section:pf-main-thm:1} as \Lemmas \ref{lemma:hegedHus2010} and \ref{lemma:R}.
We use the polynomial method, which has been widely applied in related problems, such as balancing sets of vectors and covers of the Boolean hypercube, e.g., \cite{alon1999combinatorial,alon2002balancing,alon1993covering,alon2020unbalancing,hegedHus2010balancing}.
\Lemma \ref{lemma:hegedHus2010} identifies sufficient conditions on a polynomial \(  \f \in \F_{\p}[\x_{1},\dots,\x_{4\p}]  \) 
for its degree to be lower bounded by \(  \p  \), while \Lemma \ref{lemma:R} provides a means to construct a polynomial \(  \f  \) which satisfies these conditions and relates splitting systems to the zeros of $\f$, leading to the conclusion that \(  \m \gtrsim \p \asymp \t  \).

We use \Lemma \ref{lemma:R}
to determine the prime \(  \p \asymp \t  \) and a particular collection \(  \mathcal{S}  \) of \(  \t  \)-subsets. 
In the proof of \Theorem \ref{thm:main-thm-ss-lb:1}, we construct \(  \m  \) linear functions \(  \fx{\ix}  \) whose zeros correspond with the subsets \(  \mathcal{S}_{i} \subseteq \mathcal{S}  \) of \(  \t  \)-subsets 
that intersect the \(  \ix\Th  \) blocks \(  \block_{\ix}  \) in exactly \(  \frac{\t}{2}  \) elements.
The zeros of the polynomial \(  \f = \prod_{\ix} \fx{\ix}  \) are the union of the zeros of the \(  \m  \) functions \(  \fx{\ix}  \) and correspond with the \(  \t  \)-subsets in \(  \bigcup_{i} \mathcal{S}_{i}  \). 
Observing that all \(  \t  \)-subsets in the collection \(  \mathcal{S}  \) must be zeros of \(  \f  \) if \(  (\Xset, \blockset) \) is a splitting system and applying \Lemma \ref{lemma:hegedHus2010} gives the theorem's lower bound on the size of a splitting system.



\subsection{Overview of the proof of \Lemma \ref{lemma:R}}
\label{section:pf-overview:lemma:R}

In order to overview the proof of the main technical lemma in this note, \Lemma \ref{lemma:R}, it is informally stated below.

\begin{lemma}[informal statement of \Lemma \ref{lemma:R}] 
\label{lemma:R:informmal} 
There exist a prime \(  \p \asymp \t  \) and a set \(  \Rset \in \binom{\Xset}{\t+\p}  \) such that \(  |\Rset \cap \block_{\ix}| - \frac{\t}{2} \nequiv 0 \pmod{\p} \) for every \(  \ix \in [\m]  \).
\end{lemma} 

Existence of \(  \p  \) and \(  \Rset  \) are each shown nonconstructively via the probabilistic method, where \(  \p  \) is first chosen so that a random set \(  \Rset  \), whose distribution is uniform over \(  (\t+\p)  \)-subsets of \(  \Xset  \), has the desired property with positive probability.
For any fixed \(  \p  \) falling in a particular interval, as per \Definition \ref{def:Q}, simple calculations verify that the condition that
\(  |\Rset \cap \block_{\ix}| - \frac{\t}{2} \nequiv 0 \pmod{\p} \)
can be simplified into four inequalities as the equivalent condition that
\(  |\Rset \cap \block_{\ix}| \notin \{ \MODPVALUES \} \).
Four linear maps, one per each \(  \yx \in \YVALUESSET  \), can be used to identify which \(  \ix\Th  \) blocks have
\(  \{ |\Rset \cap \block_{\ix}| = \frac{\t}{2} + \yx \p \}  \)
occurring within specified probabilities based on the value of \(  \frac{|\block_{\ix}|}{\n}  \).
The idea is then to select \(  \p  \) such that the number of \(  \ix\Th  \) blocks for which the event
\(  \{ |\Rset \cap \block_{\ix}| = \frac{\t}{2} + \yx \p \}  \)
occurs with probability at most of order \(  \frac{\nu}{\sqrt{\t}}  \) is at most of order \(  \smash{\sqrt{\frac{\t}{\nu}}}  \).
Such a prime \(  \p  \) can be shown to exist via the probabilistic method and bounds on the prime-counting function.
Upon fixing \(  \p  \) to have such a property, the application of the probabilistic method to \(  \Rset  \), together with a union bound calculated according to said property, completes the proof.


\subsection{Comparison to Alon \etal \cite{alon2020unbalancing} and Heged\"{u}s \cite{hegedHus2010balancing}}

There are substantial differences between prior works and this note with technical complications that required new analysis.
While extending the analysis of Alon \etal \cite{alon2020unbalancing} and Heged\"{u}s \cite{hegedHus2010balancing} to
\(  (\n,\t)  \)-splitting systems is relatively straightforward
when all blocks have cardinality \(  \lfloor \frac{\n}{2} \rfloor  \), 
this appears to not be the case for
general
block sizes when \(  \t < \frac{\n}{2}  \).
Although the generalization of the set system proposed by Alon \etal \cite{alon2020unbalancing} allows for arbitrary block sizes, their generalization differs from the definition of splitting systems considered here.
Indeed, their condition that \(  |\Yset \cap \block_{\ix}| = \frac{|\block_{\ix}|}{2}  \), as opposed to \(  |\Yset \cap \block_{\ix}| = \frac{|\Yset|}{2}  \), \(  |\Yset|=\t  \), changes the problem both conceptually and technically when the blocks do not have cardinality \(  \lfloor \frac{\n}{2} \rfloor  \).
The primary challenge arises in the main technical lemma, \Lemma \ref{lemma:R} and its analog in reference \cite{alon2020unbalancing}.
Recall from the proof overview in \Section \ref{section:pf-overview:lemma:R} that \Lemma \ref{lemma:R} establishes existence of both a prime \(  \p  \) and a subset of \(  \Xset  \) satisfying a certain property, where
the selection of \(  \p  \) is 
nontrivial.
In contrast, Alon \etal \cite{alon2020unbalancing} simply take \(  \p \asymp \frac{\n}{4}  \) in a fairly arbitrary manner, allocating most the proof only to showing existence of a certain subset of \(  \Xset  \).
The intuition for why we need to also carefully determine \(  \p  \) is illustrated in 
a counterexample
when the blocks are of cardinality \(  \frac{\n}{3}  \) and \(  \p = \frac{\t}{2}  \), which would have reduced the lower bound achievable by our method to \(  \m \gtrsim \sqrt{\t}  \).
Such an issue is avoided with the alternative condition that \(  |\Yset \cap \block_{\ix}| = \frac{|\block_{\ix}|}{2}  \) in the set system considered by Alon \etal \cite{alon2020unbalancing}.
\section{Proof of \Theorem \ref{thm:main-thm-ss-lb:1}}
\label{section:pf-main-thm:1}

\subsection{Preliminaries}
\label{section:pf-main-thm:1:prelim}

\subsubsection{Definitions and assumptions}
\label{section:pf-main-thm:1:prelim:defs-and-assumptions}

The following absolute constants appear in the proof of \Theorem \ref{thm:main-thm-ss-lb:1}.

\begin{definition}[constants appearing in this note] 
\label{def:constants} 
Let \(  \aconst \defeq \aconstvalue  \), \(  \bconst \defeq \bconstvalue  \), \(  \alphaconst = \alphaconstvalue  \), and \(  \cconst \defeq \frac{\alphaconst}{48 \aconst (2 \bconst + 1)^{2}}  \).
\end{definition} 

Additionally, the analysis works with a particular set of primes, whose purpose will become apparent in the proof of the main technical lemma of this work (\Lemma \ref{lemma:R}).

\begin{definition}[set of primes \(  \pset  \) appearing in this note] 
\label{def:Q} 
Denote by \(  \pset \subset \Z \cap (\minp, \maxp]  \) the set of all primes in the interval \(  (\minp, \maxp]  \):
\(  \pset \defeq \{ q \in \Z \cap (\minp, \maxp] : q \text{ is prime} \}  \).
\end{definition} 

\subsubsection{Main technical lemmas}
\label{section:pf-main-thm:1:prelim:main-technical-lemmas}

The proof of \Theorem \ref{thm:main-thm-ss-lb:1} relies on two technical lemmas, \Lemmas \ref{lemma:hegedHus2010} and \ref{lemma:R}.
The former is a result of Heged\"{u}s \cite{hegedHus2010balancing},
which Alon \etal \cite{alon2020unbalancing} and Heged\"{u}s \cite{hegedHus2010balancing} used for their lower bounds on the size of \(  (\n,\frac{\n}{2})  \)-splitting systems and on the size of a generalization of the set system proposed by Alon \etal \cite{alon2020unbalancing}.
The proof of the latter technical lemma, \nm{a result of this note}, is deferred to \Section \ref{section:pf-main-thm:1:pf-lemma:R}.

\begin{lemma}[\cite{hegedHus2010balancing}; see also {\cite[\Lemma 2.1]{alon2020unbalancing}}] 
\label{lemma:hegedHus2010} 
Let \(  \p  \) be a prime, and
let \(  \Uset  \) be a set of cardinality \(  4\p  \).
If a polynomial \(  \f \in \F_{\p}[x_{\Uset}]  \) in \(  4\p  \) variables, \(  \x_{\Uset} = (\x_{\ux})_{\ux \in \Uset}  \),
satisfies
\(  \f(1_{\Tset}) \neq 0  \) and \(  \f(1_{\Yset}) = 0  \) for at least one \(  \Tset \in \binom{\Uset}{3\p}  \) and for all \(  \Yset \in \binom{\Uset}{2\p}  \),
then \(  \f  \) has degree at least \(  \p  \).
\end{lemma} 

\begin{lemma} 
\label{lemma:R} 
Let \(  \Xset  \) be a set with cardinality \(  |\Xset| = \n  \).
Let \(  \cconst > 0  \) be an absolute constant specified in \Definition \ref{def:constants}, \nm{and let \(  \pset  \) be a set of primes defined in \Definition \ref{def:Q}}.
Let \(  \block_{1}, \dots, \block_{\m} \subseteq \Xset  \).
Suppose that \(  \m \leq \mlb  \).
Then, there exist \(  \p \in \pset  \) and \(  \Rset \in \binom{\Xset}{\t+\p}  \) such that \(  |\Rset \cap \block_{\ix}| - \frac{\t}{2} \nequiv 0 \pmod{\p} \) for every \(  \ix \in [\m]  \).
\end{lemma} 


\subsection{Proof of the main theorem}
\label{section:pf-main-thm:1:pf}

\begin{proof}{\Theorem \ref{thm:main-thm-ss-lb:1}}  
Suppose for the sake of contradiction that
\(  \m \leq \mlb  \)
and \(  (\Xset, \{ \block_{1}, \dots, \block_{\m} \})  \) is an \(  (\n,\t)  \)-splitting system.
Let \(  \p \in \pset  \) be a prime and \(  \Rset \in \binom{\Xset}{\t+\p}  \) be a subset, both chosen according to \Lemma \ref{lemma:R}.
Let
\(  \Sset \in \binom{\Xset \setminus \Rset}{\p}  \)
and
\(  \Tset \in \binom{\Rset}{3\p}  \)
be disjoint subsets, and
define
\(  \Uset \defeq \Sset \cup \Tset  \)
and
\(  \Vset \defeq \Rset \setminus \Tset  \),
where
\(  |\Uset| = |\Sset \cup \Tset| = |\Sset| + |\Tset| = \p + 3\p = 4\p  \) and
\(  |\Vset| = |\Rset \setminus \Tset| = |\Rset| - |\Tset| = (\t+\p) - 3\p = \t-2\p  \)
by design.

Define \(  \m  \) linear functions
\(  \fx{1}, \dots, \fx{\m} \in \F_{\p}[\x_{\Uset}]  \)
over \(  4\p  \) variables, \(  \x_{\Uset} = (\x_{\ux})_{\ux \in \Uset}  \), which are in one-to-one correspondence with the elements of \(  \Uset  \), where for each \(  \ix \in [\m]  \),
\begin{gather*}
  \fx{\ix}(\x) = \sum_{\ux \in \Uset \cap \block_{\ix}} \x_{\ux} + |\Vset \cap \block_{\ix}| - \frac{\t}{2}
.\end{gather*}
Then, define the polynomial
\(  \f \in \F_{\p}[\x_{\Uset}]  \)
by
\begin{gather*}
  \f = \prod_{\ix=1}^{\m} \fx{\ix}
,\end{gather*}
whose degree is at most
\(  \deg(\f) \leq \m  \).
%
%
For every \(  \ix \in [\m]  \),
\begin{align*}
  \fx{\ix}(1_{\Tset})
  &= \sum_{\ux \in \Uset \cap \block_{\ix}} (1_{\Tset})_{\ux} + |\Vset \cap \block_{\ix}| - \frac{\t}{2}
  \\
  &= |\Tset \cap \Uset \cap \block_{\ix}| + |\Vset \cap \block_{\ix}| - \frac{\t}{2}
  = |\Tset \cap \block_{\ix}| + |\Vset \cap \block_{\ix}| - \frac{\t}{2}
  = |(\Tset \cup \Vset) \cap \block_{\ix}| - \frac{\t}{2}
  = |\Rset \cap \block_{\ix}| - \frac{\t}{2}
 \neq 0
,\end{align*}
where the last line is by the design of \(  \Rset  \) in \Lemma \ref{lemma:R}.
Hence,
\(  \fx{\ix} \nequiv 0  \)
for each \(  \ix \in [\m]  \).
This implies that \(  \f  \) is also nontrivial,
and moreover, taking the product of the \(  \m  \) nonzero terms \(  \fx{\ix}(1_{\Tset})  \), it follows that
\(  \f(1_{\T}) = \prod_{\ix=1}^{\m} \fx{\ix}(1_{\Tset}) \neq 0  \).

On the other hand, by assumption, \(  (\Xset, \{ \block_{1}, \dots, \block_{\m} \})  \) is an \(  (\n,\t)  \)-splitting system.
Additionally, for any \(  \Yset \in \binom{\Uset}{2\p}  \), we have \(  |\Yset \cup \Vset| = |\Yset| + |\Vset| = 2\p + (\t-2\p) = \t  \)
since \(  \Yset \cap \Vset \subseteq \Uset \cap \Vset = \emptyset  \),
and \(  |\Vset| = \t-2\p  \) by an earlier remark.
Hence, such a set \(  \Yset \cup \Vset  \) must be split by at least one of the \(  \m  \) blocks---that is, for every \(  \Yset \in \binom{\Uset}{2\p}  \), there exists \(  \ix \in [\m]  \) such that \(  |(\Yset \cup \Vset) \cap \block_{\ix}| = \frac{\t}{2}  \).
This implies that \(  \fx{\ix}(1_{\Yset}) = 0  \) since
\begin{align*}
  \fx{\ix}(1_{\Yset})
  &= \sum_{\ux \in \Uset \cap \block_{\ix}} (1_{\Yset})_{\ux} + |\Vset \cap \block_{\ix}| - \frac{\t}{2}
  \\
  &= |\Yset \cap \Uset \cap \block_{\ix}| + |\Vset \cap \block_{\ix}| - \frac{\t}{2}
  \\
  &= |\Yset \cap \block_{\ix}| + |\Vset \cap \block_{\ix}| - \frac{\t}{2}
  = |(\Yset \cup \Vset) \cap \block_{\ix}| - \frac{\t}{2}
  = \frac{\t}{2} - \frac{\t}{2}
 = 0
.\end{align*}
It follows that \(  \f  \) vanishes on \(  1_{\Yset}  \) for all \(  \Yset \in \binom{\Xset}{2\p}  \).

Therefore, by \Lemma \ref{lemma:hegedHus2010}, the degree of \(  \f  \) is lower bounded by \(  \deg(\f) \geq \p  \).
However, this contradicts the upper bound
stated earlier: \(  \deg(\f) \leq \m < \p  \).
This completes the proof of \Theorem \ref{thm:main-thm-ss-lb:1}, assuming correctness of \Lemma \ref{lemma:R}, which is proved in the next section.
\end{proof} 

\subsection{Proof of \Lemma \ref{lemma:R}}
\label{section:pf-main-thm:1:pf-lemma:R}

\begin{proof}{\Lemma \ref{lemma:R}} 
Let \(  \p \in \pset  \) be a prime, which will be chosen in \Section \ref{section:pf-main-thm:1:pf-lemma:R:p}, and let \(  \Rset \in \binom{\Xset}{\t+\p}  \).
The following claim establishes an equivalent form of the condition which \(  \Rset  \) should avoid.

\begin{claim} 
\label{claim:eqn:pf:lemma:R:1b} 
Fix \(  \p \in \pset  \) and \(  \Rset \in \binom{\Xset}{\t+\p}  \).
For \(  \ix \in [\m]  \),
\(  |\Rset \cap \block_{\ix}| - \frac{\t}{2} \equiv 0 \pmod{\p}  \)
if and only if
\(  |\Rset \cap \block_{\ix}| \in \{ \MODPVALUES \}  \).
\end{claim} 

\begin{subproof}{\Claim \ref{claim:eqn:pf:lemma:R:1b}} 
The first equation in the claim is equivalent to
\(  |\Rset \cap \block_{\ix}| \in \{ \frac{\t}{2} + \yx \p : \yx \in \Z \}  \),
and hence, the claim will follow from showing that \(  |\Rset \cap \block_{\ix}| \neq \frac{\t}{2} + \yx \p  \) for any \(  \yx \in \Z \setminus \YVALUESSET  \).
Since \(  |\Rset| = \t+\p  \) and \(  \p > \minp  \), the inequality
\(  \frac{\t}{2}+\yx\p > \t+\p \geq |\Rset \cap \block_{\ix}|  \)
holds for all \(  \ix \in [\m]  \) if \(  \yx \geq 3  \).
In the other direction,
\(  \frac{\t}{2}+\yx\p < 0 \leq |\Rset \cap \block_{\ix}|  \)
if \(  \yx \leq -2  \).
Thus, for each \(  \ix \in [\m]  \),
\(  |\Rset \cap \block_{\ix}| = \frac{\t}{2} + \yx \p  \)
only if
\(  \yx \in (-2,3)  \),
where \(  (-2,3) \cap \Z = \YVALUESSET  \).
The claim follows.
\end{subproof} 

This suggests a strategy through the following lemma and the probabilistic method.

\begin{lemma} 
\label{lemma:lemma:R:main-ineq} 
There exists \(  \p \in \pset  \) such that
\begin{gather}
\label{eqn:pf:lemma:R:1b}
  \sum_{\ix=1}^{\m} \Pr\Bigl( \Bigl\{ |\Rset \cap \block_{\ix}| \in \Bigl\{ \MODPVALUES \Bigr\} \Bigr\} \Bigr)
  < 1
,\end{gather}
where the randomness is with respect to \(  \Rset \sim \Uniform( \binom{\Xset}{\t+\p} )  \).
\end{lemma} 

Before proving \Lemma \ref{lemma:lemma:R:main-ineq}, we show that it leads to the conclusion of the proof of \Lemma \ref{lemma:R}.
Let \(  \Rset \sim \Uniform( \binom{\Xset}{\t+\p} )  \).
By \Claim \ref{claim:eqn:pf:lemma:R:1b} and a union bound, \Equation \eqref{eqn:pf:lemma:R:1b} in \Lemma \ref{lemma:lemma:R:main-ineq} implies
\begin{gather}
\label{eqn:pf:lemma:R:1c}
  \Pr\Bigl( \bigcup_{\ix=1}^{\m} \Bigl\{ |\Rset \cap \block_{\ix}| - \frac{\t}{2} \equiv 0 \pmod{\p} \Bigr\} \Bigr)
  < 1
,\end{gather}
and by complementation,
\begin{gather}
\label{eqn:pf:lemma:R:1d}
  \Pr\Bigl( \bigcap_{\ix=1}^{\m} \Bigl\{ |\Rset \cap \block_{\ix}| - \frac{\t}{2} \nequiv 0 \pmod{\p} \Bigr\} \Bigr)
  > 0
.\end{gather}
Therefore, by the probabilistic method and \Equation \eqref{eqn:pf:lemma:R:1d}, there exists \(  \Rset \in \binom{\Xset}{\t+\p}  \) such that
\(  |\Rset \cap \block_{\ix}| - \frac{\t}{2} \nequiv 0 \pmod{\p}  \)
for all \(  \ix \in [\m]  \).
\end{proof} 

It remains to prove \Lemma \ref{lemma:lemma:R:main-ineq}.

\begin{proof}{\Lemma \ref{lemma:lemma:R:main-ineq}} 
The proof of \Lemma \ref{lemma:lemma:R:main-ineq} is divided across \Sections \ref{section:pf-main-thm:1:pf-lemma:R:notations}--\ref{section:pf-main-thm:1:pf-lemma:R:conclusion}, which break the argument down into three primary steps.
\Section \ref{section:pf-main-thm:1:pf-lemma:R:notations} defines additional notations used in the proof.
Then, \Section \ref{section:pf-main-thm:1:pf-lemma:R:p}, determines a prime \(  \p \in \pset  \) that satisfies
\nm{a key property,}
and \Section \ref{section:pf-main-thm:1:pf-lemma:R:prob-bound} upper bounds the probability of the event
\(  \{ |\Rset \cap \block_{\ix}| \in \{ \MODPVALUES \} \}  \)
for each individual \(  \ix \in [\m]  \)  given the choice of the prime \(  \p  \) with the
\nm{property}
obtained in \Section \ref{section:pf-main-thm:1:pf-lemma:R:p}.
Finally, \Section \ref{section:pf-main-thm:1:pf-lemma:R:conclusion} concludes the proof by deriving the bound in \Equation \eqref{eqn:pf:lemma:R:1b} of \Lemma \ref{lemma:lemma:R:main-ineq}.


\subsubsection{Additional notations}
\label{section:pf-main-thm:1:pf-lemma:R:notations}

For \(  \yx \in \R  \), we consider the map \(  \wtop{\yx} : \R \to \R  \) given by
\(  \wtop{\yx}(x) = \frac{(1-2x)\t}{2(x-\yx)}  \)
for \(  x \in \R  \), and its inverse
\(  \wtop{\yx}^{-1}(z) = \frac{\t+2\yx z}{2(\t+z)}  \), \(  z \in \R  \).
We denote \(  \wtop{} \defeq \wtop{0}  \).
We observe that for each \(  \yx \in \YVALUESSET  \), the function \(  \wtop{\yx}(x)  \) is continuous and monotonic over the interval \(  (0,1)  \), and its inverse \(  \wtop{\yx}^{-1}(z)  \) is continuous and monotonic over the positive real line.

Let
\(  \w_{\ix} \defeq \frac{|\block_{\ix}|}{\n}  \)
for \(  \ix \in [\m]  \), and let
\(  \Wset \defeq \{ \w_{1}, \dots, \w_{\m} \}  \).
Set \(  \h \defeq \frac{1}{5\sqrt{\t}}  \)
and
\(  \k \defeq \lceil \frac{1}{\h} \rceil \asymp \sqrt{\t}  \).
We define a partition of \(  (0,1)  \) into \(  \k  \) disjoint intervals:
\(  \Uinterval{1} \sqcup \cdots \sqcup \Uinterval{\k} = (0,1)  \) by setting
\(  \Uinterval{\jx} \defeq ((\jx-1)\h, \jx\h] \cap (0,1)  \)
for \(  \jx \in [\k]  \).

We also consider collections of sets induced by the map \(  \wtop{\yx}^{-1} \circ \wtop{}  \) by letting
\(  \Vyinterval{\yx}{\jx} \defeq \wtop{\yx}^{-1}[\wtop{}[\Uinterval{\jx}]] \cap (0,1)  \),
for \(  \jx \in [\k]  \) and \(  \yx \in \YVALUESSET  \).
We use the convention that
\(  \Vyinterval{\yx}{\jx} \defeq \emptyset  \)
for any \(  \jx \notin [\k]  \).
Note that \(  \Vyinterval{0}{\jx} = \Uinterval{\jx}  \) for \(  \jx \in [\k]  \).
Moreover, for each \(  \yx \in \YVALUESSET  \),
\(  \Vyinterval{\yx}{1} \sqcup \cdots \sqcup \Vyinterval{\yx}{\k}  \) is a partition of \(  (0,1)  \) into (possibly empty) disjoint convex sets.
This follows from the observation that for \(  \yx \in \YVALUESSET  \),
\(  \wtop{\yx}^{-1} \circ \wtop{}(x) = x + \yx -2 x \yx  \), \(  x \in \R  \),
is a 
\nm{linear map}
with
\(  \wtop{\yx}^{-1} \circ \wtop{}[(0,1)] \supseteq (0,1)  \).

For \(  \yx \in \YVALUESSET  \) and \(  \jx,\lx \in [\k]  \), the index sets are denoted by
\(  \Iyjl{\yx}{\jx}{0} \defeq \{ \ix \in [\m] : \w_{\ix} \in \Vyinterval{\yx}{\jx-1} \cup \Vyinterval{\yx}{\jx} \cup \Vyinterval{\yx}{\jx+1} \}  \)
and
\(  \Iyjl{\yx}{\jx}{\lx} \defeq \{ \ix \in [\m] : \w_{\ix} \in \Vyinterval{\yx}{\jx-\lx-1} \cup \Vyinterval{\yx}{\jx+\lx+1} \}  \).
Again, for any \(  \jx \notin [\k]  \) or \(  \lx \notin \{0,\dots,\k\}  \), set
\(  \Iyjl{\yx}{\jx}{\lx} \defeq \emptyset  \).
For \(  \yx \in \YVALUESSET  \), \(  \Iyjl{\yx}{\jx}{0} \sqcup \cdots \sqcup \Iyjl{\yx}{\jx}{\k} = [\m]  \) is a partition of \(  [\m]  \) into (possibly empty) disjoint sets due to the property stated earlier for the sets \(  \Vyinterval{\yx}{1}, \dots, \Vyinterval{\yx}{\k}  \).

\subsubsection{Selecting \(  \p  \)}
\label{section:pf-main-thm:1:pf-lemma:R:p}

This section establishes existence a prime \(  \p \in \pset  \) which satisfies a key property, specified \nm{as condition \ref{property:j0:b}} in \Lemma \ref{lemma:j0,p}, 
that can be used later to bound
the probability of the union of events
\(  \bigcup_{\ix=1}^{\m} \{ |\Rset \cap \block_{\ix}| \in \{ \MODPVALUES \} \}  \).

\begin{lemma} 
\label{lemma:j0,p} 
Let \(  \aconst, \bconst, \cconst > 0  \) be absolute constants specified in \Definition \ref{def:constants}.
Suppose \(  \m \leq \mlb  \).
Then, there exist \(  \jO \in [\k]  \) and \(  \p \in \wtop{}[\Uinterval{\jO}] \cap \pset  \) such that
the following condition
\Enum[\label{property:j0:b}]{A}
holds:
for every \(  \yx \in \YVALUESSET  \) and for every \(  \lx \in \{0,\dots,\k\}  \), \(  |\Iyjl{\yx}{\jO}{\lx}| \leq \frac{1}{\NYVALUES \aconst(2\bconst+1)} e^{\frac{\lx}{2\bconst}} \sqrt{\TMARGIN}  \).
\end{lemma} 


\begin{subproof}{\Lemma \ref{lemma:j0,p}} 
The proof consists of the following two stages.
\Enum[\label{step:pf:lemma:j0,p:i}]{i}
Denote the set of 
indexes for which \(  \wtop{}[\Uinterval{\jx}]  \) has a nonempty intersection with \(  \pset  \) by \(  \jOset \defeq \{ \jx \in [\k] : \wtop{}[\Uinterval{\jx}] \cap \pset \neq \emptyset \}  \).
We will establish a lower bound of
\(  |\jOset| \geq \alphaconst \sqrt{\TMARGIN} \asymp \sqrt{\t}  \),
where \(  \alphaconst = \alphaconstvalue  \),
using bounds for the prime-counting function.
\Enum[\label{step:pf:lemma:j0,p:ii}]{ii}
Then, through the probabilistic methods, we will prove the existence of at least one \(  \jx \in \jOset  \) satisfying \ref{property:j0:b}.
To do so, define the random variables
\(  \Zyjl{\yx}{\jx}{\lx} = \frac{\NYVALUES \aconst (2 \bconst + 1)}{\sqrt{\TMARGIN}} e^{-\frac{\lx}{2\bconst}} |\Iyjl{\yx}{\jx}{\lx}|  \)
for \(  \yx \in \YVALUESSET  \) and \(  \lx \in \{0, \dots, \k\}  \),
where \(  \jx  \) is chosen uniformly and independently from \(  \jOset  \).
Writing \(  \Zyj{\yx}{\jx} = \sum_{\yx \in \YVALUESSET} \sum_{\lx \in \{0,\dots,\k\}} \Zyjl{\yx}{\jx}{\lx}  \), it will be shown that
\nm{\(  \E[\Zyj{\yx}{\jx}] < 1  \)}
when \(  \m  \) is below the threshold stated in \Lemma \ref{lemma:j0,p}.
This will imply, by the probabilistic method and the nonnegativity of the random variables \(  \Zyjl{\yx}{\jx}{\lx}  \), \(  \yx \in \YVALUESSET  \), \(  \lx \in \{0,\dots,\k\}  \), that there exists \(  \jx \in \jOset  \) such that
\(  |\Iyjl{\yx}{\jx}{\lx}| \leq \frac{1}{\NYVALUES \aconst(2\bconst+1)} e^{\frac{\lx}{2\bconst}} \sqrt{\TMARGIN}  \)
for all \(  \yx \in \YVALUESSET  \) and all \(  \lx \in \{ 0,\dots,\k \}  \) simultaneously.

\paragraph{\Step \ref{step:pf:lemma:j0,p:i}} 

For any \(  \xv > 0  \), the prime-counting function, \(  \pi(\xv)  \), specifies the number of primes in the interval \(  (0,\xv]  \).
Recalling that \(  \pset \subseteq (\minp,\maxp]  \),
the lower bound on \(  |\jOset|  \) will follow from lower bounding the cardinality of \(  \pset  \),
and upper bounding the number of primes contained in any single \(  \jx\Th  \) interval, \(  \wtop{}[\Uinterval{\jx}]  \), \(  \jx \in [\k]  \), for which \(  \wtop{}[\Uinterval{\jx}] \cap (\minp,\maxp] \neq \emptyset  \).
We apply the following facts about prime numbers.


\begin{lemma}[\cite{dusart2018explicit}] 
\label{lemma:pi-lb-ub} 
Let \(  \xv \geq 88789  \).
Then,
\begin{gather}
\label{eqn:lemma:pi-lb-ub:1}
  \frac{\xv}{\log(\xv)}
  \Bigl( 1 + \frac{1}{\log(\xv)} + \frac{2}{\log^{2}(\xv)} \Bigr)
  \leq
  \pi(\xv)
  \leq
  \frac{\xv}{\log(\xv)}
  \Bigl( 1 + \frac{1}{\log(\xv)} + \frac{2}{\log^{2}(\xv)} + \frac{7.59}{\log^{3}(\xv)} \Bigr)
.\end{gather}
\end{lemma} 



\begin{lemma}[\cite{brun1920crible}] 
\label{lemma:pi(x+y)-p(x)} 
Let \(  \xv, \yv > 0  \), where \(  \xv \geq \xv_{0}  \) for an absolute constant \(  \xv_{0} > 0  \).
Then,
\begin{gather}
\label{eqn:lemma:pi(x+y)-p(x):1}
  \pi(\xv+\yv) - \pi(\xv) < \frac{8 \yv}{\log(\yv)}
.\end{gather}
\end{lemma} 


By \Lemma \ref{lemma:pi-lb-ub},
\begin{align*}
  |\pset|
  &=
  \pi \Bigl( \maxp \Bigr) - \pi \Bigl( \Bigl\lfloor \minp \Bigr\rfloor+1 \Bigr)
  \\
  &\geq
  \frac{5\t}{12 \log( \maxp )}
  \Bigl( 1 + \frac{1}{\log( \maxp )} + \frac{2}{\log^{2}( \maxp )} \Bigr)
  -
  \frac{\minp+1}{\log( \minp )}
  \Bigl( 1 + \frac{1}{\log( \minp )} + \frac{2}{\log^{2}( \minp )} + \frac{7.59}{\log^{3}( \minp )} \Bigr)
  \geq
  \frac{\t}{6\log(\t)}
,\end{align*}
which holds for \(  \t \geq 450  \).

Observe that
\(  \frac{d}{dz} \wtop{}(z), \frac{d}{dz} \wtop{}^{-1}(z) < 0  \) and \(  \frac{d^{2}}{dz^{2}} \wtop{}(z) > 0  \) for \(  z > 0  \).
Hence,
\(  \max \{ \sup \wtop{}[\Uinterval{\jx}] \cap (\minp,\maxp] - \inf \wtop{}[\Uinterval{\jx}] \cap (\minp,\maxp] : \jx \in [\k] \} \leq \maxp - \wtop{}(\wtop{}^{-1}(\maxp)+\h) = \frac{289\h\t}{72+204\h} \leq \sqrt{\t}  \).
Therefore, by \Lemma \ref{lemma:pi(x+y)-p(x)},
\begin{align*}
  \max \{ | \wtop{}[\Uinterval{\jx}] \cap \pset | : \jx \in [\k] \}
  &\leq
  \pi( \xv + \sqrt{\t} ) - \pi( \xv )
  <
  \frac{16 \sqrt{\t}}{\log(\t)}
\end{align*}
for some \(  \xv \geq \xv_{0}  \). 
It follows that
\begin{align*}
  |\jOset|
  &=
  |\{ \jx \in [\k] : \wtop{}[\Uinterval{\jx}] \cap \pset \neq \emptyset \}|
  \\
  &\geq
  \frac{|\pset|}
       {\max \{ | \wtop{}[\Uinterval{\jx}] \cap \pset | : \jx \in [\k] \}}
  \geq
  \frac{\t}{6\log(\t)} \frac{\log(\t)}{16 \sqrt{\t}}
  =
  \alphaconst \sqrt{\TMARGIN}
,\end{align*}
where \(  \alphaconst \defeq \alphaconstvalue  \).
This completes \Step \ref{step:pf:lemma:j0,p:i}.

\paragraph{\Step \ref{step:pf:lemma:j0,p:ii}} 

We follow the notations from \Sections \ref{section:pf-main-thm:1:pf-lemma:R:notations} and \ref{section:pf-main-thm:1:pf-lemma:R:p}.

\begin{claim} 
\label{claim:pf:lemma:j0,p:sum|I_jl^y|} 
For any \(  \yx \in \YVALUESSET  \) and \(  \lx \in \{0,\dots,\k\}  \),
\(
  \sum_{\jxx \in \jOset}
  |\Iyjl{\yx}{\jxx}{\lx}|
  \leq
  3\m
\).
\end{claim} 

\begin{subproof}{\Claim \ref{claim:pf:lemma:j0,p:sum|I_jl^y|}} 
Fix \(  \yx \in \YVALUESSET  \) arbitrarily, and suppose \(  \lx=0  \).
For \(  \jxx \in [\m]  \), let
\(  \IyjlX{\yx}{\jxx} \defeq \{ \ix \in [\m] : \w_{\ix} \in \Vyinterval{\yx}{\jxx} \}  \),
and for any \(  \jxx \notin [\k]  \), set
\(  \IyjlX{\yx}{\jxx} \defeq \emptyset  \).
The set \(  \Iyjl{\yx}{\jxx}{\lx}  \) can be partitioned into three pairwise disjoint sets as
\(  \Iyjl{\yx}{\jxx}{\lx} = \Iyjl{\yx}{\jxx}{0} = \IyjlX{\yx}{\jxx-1} \sqcup \IyjlX{\yx}{\jxx} \sqcup \IyjlX{\yx}{\jxx+1}  \)
by the definition of \(  \Iyjl{\yx}{\jxx}{0}  \) and disjointness of \(  \Vyinterval{\yx}{1}, \dots, \Vyinterval{\yx}{\k}  \).
Then,
\begin{align*}
  \sum_{\jxx \in \jOset}
  |\Iyjl{\yx}{\jxx}{\lx}|
  &\leq
  \sum_{\jxx=1}^{\k}
  |\Iyjl{\yx}{\jxx}{\lx}|
  =
  \sum_{\jxx=1}^{\k}
  ( |\IyjlX{\yx}{\jxx-1}| + |\IyjlX{\yx}{\jxx}| + |\IyjlX{\yx}{\jxx+1}| )
  \\
  &=
  \sum_{\jxx=1}^{\k}
  |\IyjlX{\yx}{\jxx-1}|
  +
  \sum_{\jxx=1}^{\k}
  |\IyjlX{\yx}{\jxx}|
  +
  \sum_{\jxx=1}^{\k}
  |\IyjlX{\yx}{\jxx+1}|
\leq
  3
  \sum_{\jxx=1}^{\k}
  |\IyjlX{\yx}{\jxx}|
  =
  3
  \sum_{\jxx=1}^{\k}
  \sum_{\ix=1}^{\m}
  \I( \ix \in \IyjlX{\yx}{\jxx} )
  \\
  &=
  3
  \sum_{\ix=1}^{\m}
  \sum_{\jxx=1}^{\k}
  \I( \ix \in \IyjlX{\yx}{\jxx} )
  =
  3
  \sum_{\ix=1}^{\m}
  \I \Bigl( \ix \in \bigcup_{\jxx=1}^{\k} \IyjlX{\yx}{\jxx} \Bigr)
  \leq
  3
  \sum_{\ix=1}^{\m}
  \I( \ix \in [\m] )
  =
  3\m
.\end{align*}
The case for \(  \lx \in [\k]  \) can be argued analogously,
but each \(  \Iyjl{\yx}{\jxx}{\lx}  \), \(  \jxx \in [\k]  \), is partitioned into two parts, 
yielding a bound of \(  \sum_{\jxx \in \jOset} |\Iyjl{\yx}{\jxx}{\lx}| \leq 2\m  \) for each \(  \lx \in [\k]  \).
Taking the maximum upper bound across all \(  \lx \in \{0,\dots,\k\}  \), the claim follows.
\end{subproof} 

\begin{claim} 
\label{claim:pf:lemma:j0,p:E[Z_j]} 
If \(  \m \leq \mlb  \), then
the expected value of \(  \Zyj{\yx}{\jx}  \) is strictly less than \(  1  \).
\end{claim} 

\begin{subproof}{\Claim \ref{claim:pf:lemma:j0,p:E[Z_j]}} 
The expected value of \(  \Zyj{\yx}{\jx}  \) 
is bounded from above by \(  1  \) as shown in the following where
the expectation is with respect to the random variable \(  \jx  \):
\begin{align*}
  \E[\Zyj{\yx}{\jx}{}]
  &=
  \frac{1}{|\jOset|}
  \sum_{\jxx \in \jOset}
  \Zyj{\yx}{\jxx}
  =
  \frac{1}{|\jOset|}
  \sum_{\yx \in \YVALUESSET}
  \sum_{\jxx \in \jOset}
  \sum_{\lx=0}^{\k}
  \Zyjl{\yx}{\jxx}{\lx}
  \\
  &\leq
  \frac{1}{\alphaconst \sqrt{\t}}
  \sum_{\yx \in \YVALUESSET}
  \sum_{\jxx \in \jOset}
  \sum_{\lx=0}^{\k}
  \frac{\NYVALUES \aconst (2 \bconst + 1)}{\sqrt{\TMARGIN}} e^{-\frac{\lx}{2\bconst}} |\Iyjl{\yx}{\jxx}{\lx}|
  \\
  &=
  \frac{\NYVALUES \aconst (2 \bconst + 1)}{\alphaconst \TMARGIN}
  \sum_{\yx \in \YVALUESSET}
  \sum_{\lx=0}^{\k}
  \frac{}{} e^{-\frac{\lx}{2\bconst}}
  \sum_{\jxx \in \jOset}
  |\Iyjl{\yx}{\jxx}{\lx}|
  \\
  &\leq
  \frac{12 \aconst (2 \bconst + 1) \m}{\alphaconst \TMARGIN}
  \sum_{\yx \in \YVALUESSET}
  \sum_{\lx=0}^{\k}
  \frac{}{} e^{-\frac{\lx}{2\bconst}}
  \\
  &=
  \frac{48 \aconst (2 \bconst + 1) \m}{\alphaconst \TMARGIN}
  \sum_{\lx=0}^{\k}
  e^{-\frac{\lx}{2\bconst}}
  <
 \frac{48 \aconst (2 \bconst + 1) \m}{\alphaconst \TMARGIN}
 \frac{1}{1-e^{-\frac{1}{2\bconst}}}
  \leq
  \frac{48 \aconst (2 \bconst + 1)^{2} \m}{\alphaconst \TMARGIN}
  \leq 1
,\end{align*}
where we apply \Claim \ref{claim:pf:lemma:j0,p:sum|I_jl^y|} and use the fact $(1+ \frac1{2\bconst} )^{-1}\ge e^{-\frac1{2\bconst}}$, and the last line holds by the assumption that
\(  \m \leq \cconst \t  \),
where
\(  \cconst = \frac{\alphaconst}{48 \aconst (2 \bconst + 1)^{2}}  \).
\end{subproof} 

By the probabilistic method, \Claim \ref{claim:pf:lemma:j0,p:E[Z_j]} implies that there exists \(  \jx \in \jOset  \) such that condition \ref{property:j0:b} of \Lemma \ref{lemma:j0,p} holds.
This completes the proof of the lemma.
\end{subproof} 

\subsubsection{Bounding the probability of the events \(  \{ |\Rset \cap \block_{\ix}| \in \{ \MODPVALUES \} \}  \)}
\label{section:pf-main-thm:1:pf-lemma:R:prob-bound}

With the selection of \(  \p \in \pset  \) from \Section \ref{section:pf-main-thm:1:pf-lemma:R:p}, which satisfies condition
\ref{property:j0:b} in \Lemma \ref{lemma:j0,p}, the probabilities of the events \(  \{ |\Rset \cap \block_{\ix}| \in \{ \MODPVALUES \} \}  \), \(  \ix \in [\m]  \), are \nm{each} bounded from above as follows.

\begin{lemma} 
\label{lemma:single-prob-ub} 
We follow the notation from \Sections \ref{section:pf-main-thm:1:pf-lemma:R:notations} and \ref{section:pf-main-thm:1:pf-lemma:R:p}.
Then, for each \(  \ix \in \Iyjl{\yx}{\jO}{\lx}  \),
\begin{gather}
  \Pr\Bigl( \Bigl\{ |\Rset \cap \block_{\ix}| = \frac{\t}{2}+\yx\p \Bigr\} \Bigr)
  \leq
  \frac{\aconst e^{-\frac{\lx}{\bconst}}}{\sqrt{\TMARGIN}}
.\end{gather}
\end{lemma} 

The proof of \Lemma \ref{lemma:single-prob-ub} depends on expressions involving binomial coefficients and on a tail bound for the hypergeometric distribution, stated below in \Lemmas \ref{lemma:binom-coef-ineq:1}--\ref{lemma:hypergeometric-tail-bound}.
The proofs of \Lemmas \ref{lemma:binom-coef-ineq:1} and \ref{lemma:binom-coef-ineq:2}  follow from Stirling's approximation.
A derivation of an equivalent form for \Lemma \ref{lemma:hypergeometric-tail-bound} can be found in \cite{skala2013hypergeometric}.


\begin{lemma} 
\label{lemma:binom-coef-ineq:1} 
Let \(  \rv,\sv \in \Z_{+}  \), where \(  \sv < \rv  \), and let \(  \uv \in [0,1]  \).
Let \(  0< \crsu{\rv}{\sv}{\uv} \leq e^{1/12}  \) be a constant. 
Then,
\begin{gather}
\label{eqn:lemma:binom-coef-ineq:1:1}
  \sup_{\vv \in [0,1]}
  \frac{\binom{\lfloor \uv \rv \rfloor}{\lfloor \vv \sv \rfloor} \binom{\lceil (1-\uv) \rv \rceil}{\lceil (1-\vv) \sv \rceil}}
       {\binom{\rv}{\sv}}
  \leq
  \frac{\crsu{\rv}{\sv}{\uv}}
       {\sqrt{2\pi \uv (1-\uv) \sv}}
  \sqrt{\frac{\rv}{\rv-\sv}}
.\end{gather}
\end{lemma} 

\begin{lemma} 
\label{lemma:binom-coef-ineq:2} 
Let \(  \rv,\sv \in \Z_{+}  \), where \(  \sv < \frac{3\rv}{4}  \), and let \(  \uv \in [0,\frac{1}{2}]  \) and \(  \vv \in [-1,1]  \), where
\(  \uv \asymp 1  \) and
\(  |\vv| \leq \frac{\uv}{3}  \).
Let \( 0< \crsuv{\rv}{\sv}{\uv}{\vv} \leq e^{1/12}  \) be a constant. 
Then,
\begin{gather}
\label{eqn:lemma:binom-coef-ineq:2:1}
  \frac{\binom{\lfloor \uv \rv \rfloor }{\lfloor (\uv+\vv) \sv \rfloor} \binom{\lceil (1-\uv) \rv \rceil}{\lceil (1-\uv-\vv) \sv \rceil}}
       {\binom{\rv}{\sv}}
  \leq
  \frac{\crsuv{\rv}{\sv}{\uv}{\vv} e^{-\frac{\vv^{2} \sv}{3}}}
       {\sqrt{2\pi \uv (1-\uv) \sv}}
  \sqrt{\frac{\rv}{\rv-\sv}}
  \frac{1}{\sqrt{(1-\frac{\vv\sv}{\uv(\rv-\sv)}) (1+\frac{\vv\sv}{(1-\uv)(\rv-\sv)})}}
.\end{gather}
\end{lemma}


\begin{lemma}[\cite{skala2013hypergeometric}] 
\label{lemma:hypergeometric-tail-bound} 
Let \(  \rv, \sv \in \Z_{+}  \), \(  \uv \in \{0, \frac{1}{\rv}, \frac{2}{\rv}, \dots, 1\}  \), and \(  \vv \in [0,1]  \), where \(  \rv \geq \sv  \).
Let \(  \Zset  \) be a set of cardinality \(  |\Zset| = \rv  \), and fix \(  \Aset \in \binom{\Zset}{\uv \rv}  \).
Let \(  \Bset \sim \Uniform( \binom{\Zset}{\sv} )  \).
Then,
\begin{gather}
\label{eqn:lemma:hypergeometric-tail-bound:1}
  \max \{
    \Pr( |\Aset \cap \Bset| - \uv \sv \leq -\vv \sv ),
    \Pr( |\Aset \cap \Bset| - \uv \sv \geq \vv \sv )
  \}
  \leq
  e^{-2 \vv^{2} \sv}
.\end{gather}
\end{lemma} 
\begin{subproof}{\Lemma \ref{lemma:single-prob-ub}} 
Arbitrarily fix \(  \yx \in \YVALUESSET  \) and \(  \lx \in \{0,\dots,\k\}  \), such that \(  \Iyjl{\yx}{\jO}{\lx} \neq \emptyset  \), and consider any \(  \ix \in \Iyjl{\yx}{\jO}{\lx}  \).
The event
\(  \{ |\Rset \cap \block_{\ix}| = \frac{\t}{2} + \yx \p \}  \)
can be rewritten as
\(  \{ |\Rset \cap \block_{\ix}| = \w(\t+\p) \}  \)
where \(  \w = \wtop{\yx}^{-1}(\p) = \frac{\t+2\yx\p}{2(\t+\p)} \in \Vyinterval{\yx}{\jO}  \).
Note that \(  \w \geq \frac{1}{17}  \), which coincides with \(  \p = \maxp  \) and \(  \yx=-1  \).

We consider the case that \(  \lx = 0  \).
Then,
\begin{align*}
  \Pr\Bigl( \Bigl\{ |\Rset \cap \block_{\ix}| = \frac{\t}{2} + \yx \p \Bigr\} \Bigr)
  &=
  \Pr( \{ |\Rset \cap \block_{\ix}| = \w(\t+\p) \} )
  =
  \frac{\binom{\w_{\ix} \n}{\w(\t+\p)} \binom{(1-\w_{\ix}) \n}{(1-\w)(\t+\p)}}
       {\binom{\n}{\t+\p}}
  \\
  &\leq
  \frac{\crsu{\n}{\t}{\w_{\ix}}}
       {\sqrt{2\pi \w_{\ix} (1-\w_{\ix}) (\t+\p)}}
  \sqrt{\frac{\n}{\n-\t-\p}}
  \leq
  \frac{\aconst}{\sqrt{\TMARGIN}}
  =
  \frac{\aconst e^{-\frac{\ell}{\bconst}}}{\sqrt{\TMARGIN}}
\TagEqn \label{eqn:lemma:single-prob-ub:1}
,\end{align*}
where the first inequality is by \Equation \eqref{eqn:lemma:binom-coef-ineq:1:1} in \Lemma \ref{lemma:binom-coef-ineq:1}, the second inequality holds since \(  \w_{\ix} \geq \w-\h \geq \frac{1}{17} - \frac{1}{5\sqrt{\t}}  \) and \(  \aconst = \aconstvalue  \), and the last equality uses that \(  \lx = 0  \).


For the case of \(  \lx \geq 1  \), we observe that
\(  \w_{\ix}(\t+\p) - (\frac{\t}{2} + \yx \p) = \w_{\ix}(\t+\p) - \w(\t+\p) = (\w_{\ix}-\w)(\t+\p)  \).
Hence,
\(  \frac{\t}{2} + \yx \p = \w_{\ix}(\t+\p) - (\w_{\ix}-\w)(\t+\p)  \).
If \(  \w_{\ix} \leq \w  \), then
\(  \w_{\ix} \in \Vyinterval{\yx}{\jO-\lx-1}  \), and so
\(  \w-\w_{\ix} \geq \inf \Vyinterval{\yx}{\jO} - \sup \Vyinterval{\yx}{\jO-\lx-1} \geq ((\jO-1) - (\jO-\lx-1) - 1) \h = \lx \h  \).
Suppose \(  \lx \leq \sqrt{\frac{\log(\sqrt{\TMARGIN}/\aconst)}{2\h^{2}(\t+\p)}+\frac{1}{144}} + \frac{1}{12}  \).
Then, noting that \(  \t+\p < \frac{3\t}{2} \leq \frac{3\n}{4}  \), by \Lemma \ref{lemma:binom-coef-ineq:2},
\begin{align*}
  \Pr\Bigl( \Bigl\{ |\Rset \cap \block_{\ix}| = \frac{\t}{2} + \yx \p \Bigr\} \Bigr)
  &=
  \Pr( \{ |\Rset \cap \block_{\ix}| = \w(\t+\p) \} )
  =
  \Pr( \{ |\Rset \cap \block_{\ix}| = (\w_{\ix}+(\w-\w_{\ix}))(\t+\p) \} )
  \\
  &=
  \frac{\binom{\w_{\ix} \n}{(\w_{\ix}+(\w-\w_{\ix}))(\t+\p)} \binom{(1-\w_{\ix}) \n}{(1-(\w_{\ix}+(\w-\w_{\ix})))(\t+\p)}}
       {\binom{\n}{\t+\p}}
  \\
  &\leq
  \frac{\crsuv{\n}{\t+\p}{\w_{\ix}}{\lx\h} e^{-\frac{\lx^{2} \h^{2} (\t+\p)}{3}}}
       {\sqrt{2\pi \w_{\ix} (1-\w_{\ix}) (\t+\p)}}
  \sqrt{\frac{\n}{\n-\t-\p}}
  \frac{1}{\sqrt{(1-\frac{\lx \h (\t+\p)}{\w_{\ix} (\n-\t-\p)}) (1+\frac{\lx \h (\t+\p)}{(1-\w_{\ix})(\n-\t-\p)})}}
  \\
  &\leq
  \frac{\aconst e^{-\frac{\lx}{\bconst}}}{\sqrt{\TMARGIN}}
\TagEqn \label{eqn:lemma:single-prob-ub:2}
,\end{align*}
where the last inequality uses facts including \(  \aconst = \aconstvalue  \), \(  \lx \h \lesssim \sqrt{\frac{\log(\t)}{\t}}  \), and \(  \w_{\ix} \geq \w-\lx\h \gtrsim \frac{1}{17} - \sqrt{\frac{\log(\t)}{\t}}  \).


Otherwise, for \(  \lx > \sqrt{\frac{\log(\sqrt{\TMARGIN}/\aconst)}{2\h^{2}(\t+\p)}+\frac{1}{144}} + \frac{1}{12}  \), observe that
\begin{align*}
  \Pr\Bigl( \Bigl\{ |\Rset \cap \block_{\ix}| = \frac{\t}{2} + \yx \p \Bigr\} \Bigr)
  &=
  \Pr( \{ |\Rset \cap \block_{\ix}| = \w(\t+\p) \} )
  \\
  &=
  \Pr( \{ |\Rset \cap \block_{\ix}| = (\w_{\ix}+(\w-\w_{\ix}))(\t+\p) \} )
  \\
  &=
  \Pr( \{ |\Rset \cap \block_{\ix}| - \w_{\ix}(\t+\p) = (\w-\w_{\ix})(\t+\p) \} )
  \\
  &\leq
  \Pr( \{ |\Rset \cap \block_{\ix}| - \w_{\ix}(\t+\p) \geq (\w-\w_{\ix})(\t+\p) \} )
  \\
  &\leq
  \Pr( \{ |\Rset \cap \block_{\ix}| - \w_{\ix}(\t+\p) \geq \lx \h (\t+\p) \} )
.\end{align*}
Then, by \Lemma \ref{lemma:hypergeometric-tail-bound} and the assumption on \(  \lx  \),
\begin{align*}
  \Pr\Bigl( \Bigl\{ |\Rset \cap \block_{\ix}| = \frac{\t}{2} + \yx \p \Bigr\} \Bigr)
  &\leq
  \Pr( \{ |\Rset \cap \block_{\ix}| - \w_{\ix}(\t+\p) \geq \lx \h (\t+\p) \} )
  \\
  &\leq
  e^{-2 \lx^{2} \h^{2} (\t+\p)}
  <
  e^{-\frac{\lx \h^{2} (\t+\p)}{3} - \log(\sqrt{\TMARGIN}/\aconst)}
  =
  \frac{\aconst e^{-\frac{\lx \h^{2} (\t+\p)}{3}}}{\sqrt{\TMARGIN}}
  \leq
  \frac{\aconst e^{-\frac{\lx}{\bconst}}}{\sqrt{\TMARGIN}}
\TagEqn \label{eqn:lemma:single-prob-ub:3}
.\end{align*}

For \(  \w_{\ix} > \w  \), we can similarly derive the same upper bound for \(  \{ |\Rset \cap \block_{\ix}| = \frac{\t}{2} + \yx \p \}  \).

Combining all above, we have for all \(  \yx \in \YVALUESSET  \), \(  \lx \in \{0,\dots,\k\}  \), and \(  \ix \in \Iyjl{\yx}{\jO}{\lx}  \),
\begin{gather*}
  \Pr\Bigl( \Bigl\{ |\Rset \cap \block_{\ix}| = \frac{\t}{2}+\yx\p \Bigr\} \Bigr)
  \leq
  \frac{\aconst e^{-\frac{\lx}{\bconst}}}{\sqrt{\TMARGIN}}
.\end{gather*}
\end{subproof} 


\subsubsection{Deriving \Equations \eqref{eqn:pf:lemma:R:1b} to complete the proof}
\label{section:pf-main-thm:1:pf-lemma:R:conclusion}

The strict inequality in \Equation \eqref{eqn:pf:lemma:R:1b} now follows:
\begin{align*}
  &\NegAlnSp
  \sum_{\ix=1}^{\m}
  \Pr\Bigl( \Bigl\{ |\Rset \cap \block_{\ix}| \in \Bigl\{ \MODPVALUES \Bigr\} \Bigr\} \Bigr)
  \\
  &\leq
  \sum_{\yx \in \YVALUESSET}
  \sum_{\ix=1}^{\m}
  \Pr\Bigl( \Bigl\{ |\Rset \cap \block_{\ix}| = \frac{\t}{2}+\yx\p \Bigr\} \Bigr)
  =
  \sum_{\yx \in \YVALUESSET}
  \sum_{\lx=0}^{\kx} \sum_{\ix \in \Iyjl{\yx}{\jO}{\lx}}
  \Pr\Bigl( \Bigl\{ |\Rset \cap \block_{\ix}| = \frac{\t}{2}+\yx\p \Bigr\} \Bigr)
  \\
  &\stackrel{(\mathrm{a})}{\leq}
  \sum_{\yx \in \YVALUESSET}
  \sum_{\lx=0}^{\kx} \sum_{\ix \in \Iyjl{\yx}{\jO}{\lx}}
  \frac{\aconst e^{-\frac{\lx}{\bconst}}}{\sqrt{\TMARGIN}}
  =
  \sum_{\yx \in \YVALUESSET}
  \sum_{\lx=0}^{\kx}
  |\Iyjl{\yx}{\jO}{\lx}|
  \frac{\aconst e^{-\frac{\lx}{\bconst}}}{\sqrt{\TMARGIN}}
  \stackrel{(\mathrm{b})}{\leq}
  \sum_{\yx \in \YVALUESSET}
  \sum_{\lx=0}^{\kx}
  \frac{e^{\frac{\lx}{2\bconst}} \sqrt{\TMARGIN}}{\NYVALUES \aconst(2\bconst+1)}
  \frac{\aconst e^{-\frac{\lx}{\bconst}}}{\sqrt{\TMARGIN}}
  \\
  &=
  \sum_{\yx \in \YVALUESSET}
  \sum_{\lx=0}^{\kx}
  \frac{e^{-\frac{\lx}{2\bconst}}}{\NYVALUES(2\bconst+1)}
  =
  \frac{1}{2\bconst+1}
  \sum_{\lx=0}^{\kx}
  e^{-\frac{\lx}{2\bconst}}
  <
  \frac{e^{\frac{1}{2\bconst}}}{(2\bconst+1)(e^{\frac{1}{2\bconst}}-1)}
  \stackrel{(\mathrm{c})}{<} 1
,\end{align*}
where (a) follows from \Lemma \ref{lemma:single-prob-ub}, (b) follows from condition \ref{property:j0:b} of \Lemma \ref{lemma:j0,p} for the choice of \(  \jO  \), and (c) applies the inequality \(  \tfrac{e^{\frac{1}{x}}}{e^{\frac{1}{x}}-1} < x+1  \) for \(  x > 0  \) and arguments in the proof of \Claim \ref{claim:pf:lemma:j0,p:E[Z_j]}.
\end{proof} 

\section{Conclusion}

This note established that for \(  \tO \leq \t \leq \frac{\n}{2}  \), the size of an \(  (\n,\t)  \)-splitting system is necessarily at least \(  \m \gtrsim \t  \), which, to the best of our knowledge, supersedes lower bounds for general
\(  \t  \)
obtained previously.
While this is
an exponential improvement 
over prior lower bounds
\nm{and an order-\(  \sqrt{\t}  \) improvement over the information lower bound,}
there is still a
\(  (\sqrt{\t} \log \n)  \)-gap
between our lower bound, \(  \m \gtrsim \t  \), and the best known upper bound,
\nm{\(  \m \lesssim \t^{3/2} \log \n  \).}
As we conjecture that the lower bound is order-wise tight,
it would be interesting to explore whether techniques and constructions beyond those proposed so far for this problem can
improve the upper bound on the size of splitting systems for \(  \t \gg 1  \). 

Splitting systems have potential applications in  problems not previously identified in the literature.
Examples for such applications include
universal reconstruction in
threshold group testing and exact recovery of binary vectors in 1-bit compressed sensing. 
Splitting systems can be related to a necessary property of measurement designs in these problems that appears to be a bottleneck for efficient signal reconstruction. 
If the optimal size of splitting systems is the conjectured \(  \m \asymp \t  \), 
a construction for splitting systems with \(  \m \lesssim \t  \) could guide the development of optimal designs and algorithms in these applications.
This would close a similar \(  \sqrt{\t}  \)-gap between the known upper bound, \(  \m \lesssim \k^{3/2} \log \frac{\n}{\k}  \) \citep{bansal2022universal,de2020subquadratic}, and lower bound, \(  \m \gtrsim \k \log \frac{\n}{\k}  \) \citep{aldridge2019group,bansal2022universal,de2020subquadratic}.
\nm{This exemplifies the merit of splitting systems beyond theoretical interest in applications in a variety of adjacent areas.
Their utility spans problems in fields such as cryptography, statistical inference, and signal processing, including signal compression and other high-dimensional problems with sparsity.}

\bibliography{refs}

\end{document}